\documentclass{amsart}

\usepackage{fancyhdr}
\usepackage{lastpage}
\usepackage{stmaryrd,yhmath}

\newcommand{\bdis}{\begin{displaymath}}
\newcommand{\edis}{\end{displaymath}}
\newcommand{\be}{\begin{equation}}
\newcommand{\ee}{\end{equation}}
\newcommand{\mbb}{\mathbb}
\newcommand{\mcal}{\mathcal}

\newcommand{\vp}{\varphi}

\newcommand{\vth}{\vartheta}

\newcommand{\zf}{\zeta\left(\frac{1}{2}+it\right)}

\theoremstyle{definition}

\newtheorem{cor}[]{Corollary}

\theoremstyle{remark}
\newtheorem{remark}[]{Remark}

\newtheorem*{mydef1}{{\bf Theorem}}

\newtheorem*{mydef6}{{\bf Example}}

\numberwithin{equation}{section}

\begin{document}

\title{Jacob's ladders, their iterations and the new class of integrals connected with parts of the Hardy-Littlewood integral of the function $|\zf|^2$}

\author{Jan Moser}

\address{Department of Mathematical Analysis and Numerical Mathematics, Comenius University, Mlynska Dolina M105, 842 48 Bratislava, SLOVAKIA}

\email{jan.mozer@fmph.uniba.sk}


\begin{abstract}
In this paper we introduce the iterations of the Jacob's ladder and the new type of integral containing certain product of the factors $|\zeta|^2$ corresponding to the components of some disconnected
set of the critical line. Next, we obtain an asymptotic formula for this integral, its factorization and, for example, the essential generalization of two Selberg's formulae (1946). \\

Dedicated to the 500th anniversary of rabbi L\" ow.
\end{abstract}
\maketitle

\section{Introduction}

\subsection{}

Let
\be \label{1.1}
Z(t)=e^{i\vth(t)}\zf,\ \vth(t)=-\frac t2\ln\pi+\text{Im}\ln\Gamma\left(\frac 14+i\frac t2\right)
\ee
be the signal generated by the Riemann zeta-function on the critical line. Hardy and Littlewood started to study the following integral in 1918
\be \label{1.2}
\int_0^T\left|\zf\right|^2{\rm d}t=\int_0^TZ^2(t){\rm d}t ,
\ee
and they have derived the following formula (see \cite{1}, pp. 151-156)
\be \label{1.3}
\int_0^T Z^2(t){\rm d}t\sim T\ln T,\ T\to\infty .
\ee
In this direction, the Titchmarsh-Kober-Atkinson (TKA) formula
\be \label{1.4}
\int_0^\infty Z^2(t)e^{-2\delta t}{\rm d}t=\frac{c-\ln(4\pi\delta)}{2\sin\delta}+\sum_{n=0}^N c_n\delta^n+\mcal{O}(\delta^{N+\epsilon})
\ee
(see \cite{6}, p. 141) where $c$ is the Euler constant, remained as an isolated result for the period of 56 years. We have obtained in our paper \cite{3} the nonlinear integral equation
\be \label{1.5}
\int_0^{\mu[x(T)]}Z^2(t)e^{-\frac{2}{x(T)}t}{\rm d}t=\int_0^T Z^2(t){\rm d}t ,
\ee
where
\bdis
\mu(y)\geq 7y\ln y .
\edis
Each function $\mu(y)$ generates a solution
\bdis
y=\vp_\mu(T)=\vp(T) .
\edis
Namely, we have shown in \cite{3} that the following infinite set of the almost exact expressions of the Hardy-Littlewood integral (\ref{1.1}) takes place
\be \label{1.6}
\begin{split}
& \int_0^T Z^2(t){\rm d}t=\vp_1(T)\ln\vp_1(T)+(c-\ln 2\pi)\vp_1(T)+c_0+\mcal{O}\left(\frac{\ln T}{T}\right) , \\
& \vp_1(T)=\frac 12\vp(T) ,
\end{split}
\ee
where $\vp(T)$ is the Jacob's ladder, i. e. the solution of the integral equation (\ref{1.5}).

\begin{remark}
Hence, we have proved that except the asymptotic formula (\ref{1.3}) possessing an unbounded error term (it is clear that the corresponding formulae of Ingham, Titchmarsh and Balasubramanian also posses
the unbounded errors, comp. \cite{3}, Remark 1) there is an infinite set of almost exact representations (\ref{1.6}) of the Hardy-Littlewood integral (\ref{1.2}).
\end{remark}

\begin{remark}
The result (\ref{1.6}) can be formulated as follows: the Jacob's ladders $\vp_1(t)$ (comp. (\ref{1.6}) and the extension in \cite{3}, p. 415; $G[\vp(T)]$) are the asymptotic solutions of the following
transcendental equation
\bdis
\int_0^T Z^2(t){\rm d}t=V(T)\ln V(T)+(c-\ln 2\pi)V(T)+c_0 .
\edis
\end{remark}

\subsection{}

Next, we have proved (see \cite{4}, (2.5)) the asymptotic formula
\be \label{1.7}
\begin{split}
& \int_T^{T+U} \left|\zf\right|^2{\rm d}t\sim U\tan[\alpha(T,U)]\ln T , \\
& U\in\left(\left. 0,\frac{T}{\ln T}\right.\right],\ T\to\infty ,
\end{split}
\ee
where $\alpha(T,U)$ is the angle of the chord of the curve
\bdis
y=\vp_1(T)
\edis
that binds the points $[T,\vp_1(T)]$ and $[T+U,\vp_1(T+U)]$ and, in addition to this, the asymptotic formula (see \cite{4}, (8.3))
\be \label{1.8}
\begin{split}
& \int_T^{T+U_1}\left|\zeta\left( \frac 12+i\vp_1(t)\right)\right|^4\left|\zf\right|^2{\rm d}t\sim \frac{1}{2\pi^2}U_1\ln^5T, \\
& U_1=T^{7/8+2\epsilon},\ T\to\infty ,
\end{split}
\ee
(small improvements to of the exponent $\frac 78$ and, similarly, of the analogous exponents $\frac 13,\ \frac 12$ are irrelevant for our purpose.
On the contrary, the case $U\in (0,1)$ is relevant).

\begin{remark}
The formula (\ref{1.8}) is the first integral asymptotic formula of the sixth order in the theory of the function $\zf$.
\end{remark}

Let us remind that the proof of the formula (\ref{1.8}) is simultaneously the proof of the following theorem (see \cite{4}, Theorem 3, p. 219): every Jacob's ladder
\bdis
\vp_1(t)=\frac 12\vp(t)
\edis
where $\vp(t)$ is an exact solution of the nonlinear integral equation (\ref{1.5}) is the asymptotic solution of the following nonlinear integral equation
\be \label{1.9}
\int_T^{T+U}\left| \zeta\left(\frac 12+ix(t)\right)\right|^4\left|\zf\right|^2{\rm d}t=\frac{1}{2\pi^2}U_1\ln^5T .
\ee

\subsection{}

Certain motivation for the next step is the well-known multiplicative formula
\be \label{1.10}
M\left(\prod_{k=1}^n X_k\right)=\prod_{k=1}^n M(X_k)
\ee
from the theory of probability. In the formula (\ref{1.10}) $X_k,\ k=1,\dots ,n$ stand for independent random variables and $M$ stands for the population mean. \\

Let
\be \label{1.11}
\begin{split}
& y=\frac 12\vp(t)=\vp_1(t),\ \vp_1^0(t)=t,\ \vp_1^1(t)=\vp_1(t), \\
& \vp_1^2(t)=\vp_1[\vp_1(t)],\dots,\vp_1^k(t)=\vp_1[\vp_1^{k-1}(t)],\dots,t\in [T,T+U] ,
\end{split}
\ee
where $\vp_1^k(t)$ denotes the $k$-th iteration of the Jacob's ladder $\vp_1(t),\ t\geq T_0[\vp_1]$. Let us remind that the functions $\vp_1^k(t),\ k=2,\dots$ are increasing since the function
$\vp_1(t)$ is increasing. \\

In this paper we obtain an asymptotic formula for a new kind of the transcendental integrals
\be \label{1.12}
\int_T^{T+U} F[\vp_1^{n+1}(t)]\prod_{k=0}^n\left|\zeta\left(\frac 12+i\vp_1^k(t)\right)\right|^2{\rm d}t,\ U\in\left(\left. 0,\frac{T}{\ln T}\right.\right]
\ee
for every fixed $n\in \mbb{N}_0$ and for every Lebesgue-integrable function
\bdis
F(t),\ t\in \left[\vp_1^{n+1}(T),\vp_1^{n+1}(T+U)\right],\ F(t)\geq 0 \ (\leq 0) .
\edis

\begin{remark}
The integral (\ref{1.12}) is a natural multiplicative generalization of the part of the Hardy-Littlewood integral (\ref{1.2}):
\be \label{1.13}
\begin{split}
& \int_T^{T+U}\left|\zf\right|^2{\rm d}t \ \rightarrow \\
& \rightarrow \ \int_T^{T+U}\left|\zf\right|^2F\left[\vp_1^{n+1}(t)\right]\prod_{k=1}^n\left|\zeta\left(\frac 12+i\vp_1^k(t)\right)\right|^2{\rm d}t .
\end{split}
\ee
\end{remark}

The following is connected with the formula (\ref{1.12}):
\begin{itemize}
\item[(A)] if $T\to\infty$ then the disconnected set
\bdis
\bigcup_{k=0}^{n+1} \left[\vp_1^k(T),\vp_1^k(T+U) \right]
\edis
looks like a one dimensional Friedmann-Hubble expanding universe,

\item[(B)] new class of integro-iterative equations,

\item[(C)] the set of $p(n+1)-1$ factorization formulae for the integral (\ref{1.12}) generated by all proper Euler's partitions of the number $n+1$ where, for example,
\bdis
p(200)-1=3\ 972\ 999\ 029\ 387 ,
\edis
and in the general case
\bdis
p(n+1)\sim \frac{1}{4(n+1)\sqrt{3}}e^{K\sqrt{n+1}},\ K=\pi\sqrt{\frac 23},\ n\to\infty
\edis
by the well-known Hardy-Ramanujan-Rademacher formula (see \cite{2}, pp. 164, 166),

\item[(D)] an essential generalization of two Selberg's formulae.
\end{itemize}

\begin{remark}
The integral (\ref{1.12}) expresses the energy of the complicated signal
\bdis
\sqrt{\left| F\left[\vp_1^{n+1}(t)\right]\right|}\prod_{k=0}^n\left|\zeta\left(\frac 12+i\vp_1^k(t)\right)\right|,\ t\in[T,T+U],
\edis
and formulae corresponding to (A) -- (D) give some properties of this energy.
\end{remark}

\begin{remark}
The integral (\ref{1.12}) is not accessible neither for the classical methods of Hardy-Littlewood and Selberg, nor for the current methods in the theory of the Riemann zeta-function.
\end{remark}

\section{The Theorem and a new class of integro-iterative equations}

\subsection{}

The following Theorem holds true.

\begin{mydef1}
For every fixed $n\in\mbb{N}_0$ and for every Lebesgue-integrable function
\bdis
\begin{split}
& F(t),\ t\in \left[\vp_1^{n+1}(T),\vp_1^{n+1}(T+U)\right],\ F(t)\geq 0 \ (\leq 0) , \\
& \int_{\vp_1^{n+1}(T)}^{\vp_1^{n+1}(T+U)}F(t){\rm d}t\not=0
\end{split}
\edis
we have
\be \label{2.1}
\begin{split}
& \int_T^{T+U}F\left[\vp_1^{n+1}(t)\right]\prod_{k=0}^n\left|\zeta\left(\frac 12+i\vp_1^k(t)\right)\right|^2{\rm d}t\sim \\
& \sim \left\{\int_{\vp_1^{n+1}(T)}^{\vp_1^{n+1}(T+U)}F(t){\rm d}t\right\}\ln^{n+1}T,\ U\in\left(\left. 0,\frac{T}{\ln^2T}\right.\right],\ T\to\infty ,
\end{split}
\ee
and the set
\be \label{2.2}
\bigcup_{k=0}^{n+1} \left[\vp_1^k(T),\vp_1^k(T+U) \right]
\ee
has the following properties:
\be \label{2.3}
t\sim\vp_1^k(t),\ \vp_1^k(T)\geq (1-\epsilon)T,\ k=0,1,\dots,n+1,
\ee
\be \label{2.4}
\vp_1^k(T+U)-\vp_1^k(T)<\frac{1}{2n+5}\frac{T}{\ln T},\ k=1,\dots,n+1,
\ee
\be \label{2.5}
\vp_1^k(T)-\vp_1^{k+1}(T+U)>0.18\times \frac{T}{\ln T},\ k=0,1,\dots,n.
\ee
Next, in the macroscopic domain, i. e. for
\be \label{2.6}
U\in \left[ T^{1/3+\epsilon},\frac{T}{\ln^2T}\right]
\ee
we have more detailed information about the set (\ref{2.2}):
\be \label{2.7}
\left|\left[\vp_1^k(T),\vp_1^k(T+U)\right]\right|=\vp_1^k(T+U)-\vp_1^k(T)\sim U,\ k=1,\dots,n+1,
\ee
\be \label{2.8}
\vp_1^k(T)-\vp_1^{k+1}(T+U)\sim (1-c)\frac{T}{\ln T},\ k=0,1,\dots,n,
\ee
\be \label{2.9}
\begin{split}
& \rho\left\{\left[\vp_1^{k-1}(T),\vp_1^{k-1}(T+U)\right];\left[\vp_1^k(T),\vp_1^k(T+U)\right]\right\}\sim \\
& \sim (1-c)\frac{T}{\ln T},\ k=1,\dots,n+1
\end{split}
\ee
where $\rho$ denotes the distance of corresponding segments.
\end{mydef1}

\begin{remark}
First of all we have that the set (\ref{2.2}) is disconnected (see (\ref{2.5}) for every admissible $U$ (see (\ref{2.1}). The components of connectedness of the set (\ref{2.2}) are distributed from the
right to the left (see (\ref{2.5}).
\end{remark}

\begin{remark}
Bellow listed properties of the set (\ref{2.2}) hold true in the macroscopic domain (\ref{2.6}):
\begin{itemize}
\item[(a)] the components of connectedness of the set (\ref{2.2}) have asymptotically equal measures (see (\ref{2.7}), i. e. the transformations
\bdis
\vp_1^k:\ [T,T+U]\to \left[\vp_1^k(T),\vp_1^k(T+U)\right]
\edis
asymptotically preserve the measure of the segment $[T,T+U]$,

\item[(b)] the adjacent intervals of this set have asymptotically equal measures (see (\ref{2.8})
\bdis
\left|\left(\vp_1^{k+1}(T+U),\vp_1^k(T)\right)\right|\sim (1-c)\frac{T}{\ln T}.
\edis
\end{itemize}
\end{remark}

Hence, by (a) and (b) the components of connectedness of the the disconnected set (\ref{2.2}) are distributed with remarkably asymptotic regularity.

\begin{remark}
The asymptotic behavior of the disconnected set (\ref{2.2}) is as follows: if $T\to\infty$ then the components of connectedness of this set receding unboundedly each from other (see (\ref{2.5}), (\ref{2.8}))
and all together are receding to infinity. Hence, if $T\to\infty$ then the set (\ref{2.2}) behaves as a one-dimensional Friedmann-Hubble expanding universe (comp. Introduction, (A)).
\end{remark}

\subsection{}

A new class of nonlinear equations is connected wit our Theorem. Namely, we obtain from our Theorem the following corollary.

\begin{cor}
Every Jacob's ladder
\bdis
\vp_1(t)=\frac 12\vp(t)
\edis
where $\vp(t)$ is an exact solution of the nonlinear integral equation (\ref{1.5}) is an asymptotic solution of the following integro-iterative equation (comp. (\ref{1.9}))
\be \label{2.10}
\begin{split}
& \frac 1U\int_T^{T+U} F[x^{n+1}(t)]\prod_{k=0}^n\left|\zeta\left(\frac 12+ix^k(t)\right)\right|^2{\rm d}t= \\
& =\left\{\int_{x^{n+1}(T)}^{x^{n+1}(T+U)}F(t){\rm d}t\right\}\ln^{n+1}T
\end{split}
\ee
in the sense that
\bdis
\frac{\frac{1}{U}\int_T^{T+U}F[\vp_1^{n+1}(t)]\prod_{k=0}^n\left|\zeta\left(\frac 12+i\vp_1^k(t)\right)\right|^2{\rm d}t}{\int_{\vp_1^{n+1}(T)}^{\vp_1^{n+1}(T+U)}F(t){\rm d}t}\sim\ln^{n+1}T,\ T\to\infty ,
\edis
where
\bdis
x^0(t)=t,\ x^1(t)=x(t),\ x^2(t)=x(x(t)),\dots ,
\edis
i. e. the function $x^k(t)$ is the $k$-th iteration of the function $x(t)$.
\end{cor}

\begin{remark}
There are fixed-point methods and other methods of the functional analysis used to study nonlinear equations. What can be obtained by using these methods in the case of the nonlinear
integro-iterative equation (\ref{2.10}) (at least in the case $F(t)=1$)?
\end{remark}

\section{On the set of factorizations of the integral (\ref{1.12}) generated by the set of all proper partitions of $n+1$}

\subsection{}

Since (see (\ref{2.1}), $F(t)=1$)
\be \label{3.1}
\int_T^{T+U} \prod_{k=0}^n\left|\zeta\left(\frac 12+i\vp_1^k(t)\right)\right|^2{\rm d}t\sim\left\{\vp_1^{n+1}(T+U)-\vp_1^{n+1}(T)\right\}\ln^{n+1}T
\ee
then for every proper partition ($n+1=n+1$ is excluded)
\be \label{3.2}
n+1=a_{j_1}+a_{j_2}+\dots+a_{j_s},\ a_{j_l}\in [1,n],\ l=1,\dots,s
\ee
we have
\be \label{3.3}
\begin{split}
& \ln^{a_{j_l}}T\sim \\
& \sim\frac{1}{\vp_1^{a_{j_l}}(T+U)-\vp_1^{a_{j_l}}(T)}\int_T^{T+U}\prod_{k=0}^{a_{j_l}-1}\left|\zeta\left(\frac 12+i\vp_1^k(t)\right)\right|^2{\rm d}t .
\end{split}
\ee
Next, we obtain from Theorem by (\ref{3.2}), (\ref{3.3}) the following formula
\be \label{3.4}
\begin{split}
& \int_T^{T+U}F[\vp_1^{n+1}(t)]\prod_{k=0}^n\left|\zeta\left(\frac 12+i\vp_1^k(t)\right)\right|^2{\rm d}t\sim \\
& \sim\int_{\vp_1^{n+1}(T)}^{\vp_1^{n+1}(T+U)}F(t){\rm d}t\ \times \\
& \times \prod_{l=1}^s\frac{1}{\vp_1^{a_{j_l}}(T+U)-\vp_1^{a_{j_l}}(T)}\int_T^{T+U}\prod_{k=0}^{a_{j_l}-1}\left|\zeta\left(\frac 12+i\vp_1^k(t)\right)\right|^2{\rm d}t .
\end{split}
\ee
Thus, if we use the weight factors
\be \label{3.5}
\begin{split}
& g_l=\frac{U}{\vp_1^{a_{j_l}}(T+U)-\vp_1^{a_{j_l}}(T)},\ l=1,\dots,s , \\
& g_{n+1}=\frac{U}{\vp_1^{n+1}(T+U)-\vp_1^{n+1}(T)} ,
\end{split}
\ee
then we obtain from Theorem by (\ref{3.4}), (\ref{3.5}) the following corollary.

\begin{cor}
The set of all proper partitions (\ref{3.2}) generates for every fixed $n\in\mbb{N}_0$ the set of $p(n+1)-1$ factorization formulae for the weighted mean-value of the integral (\ref{1.12})
\be \label{3.6}
\begin{split}
& g_{n+1}\frac{1}{U}\int_T^{T+U} F[\vp_1^{n+1}(t)]\prod_{k=0}^n\left|\zeta\left(\frac 12+i\vp_1^k(t)\right)\right|^2{\rm d}t\sim \\
& \sim \frac{1}{\vp_1^{n+1}(T+U)-\vp_1^{n+1}(T)}\int_{\vp_1^{n+1}(T)}^{\vp_1^{n+1}(T+U)}F(t){\rm d}t\ \times \\
& \prod_{l=1}^s g_l\frac{1}{U}\int_T^{T+U}\prod_{k=0}^{a_{j_l}-1}\left|\zeta\left(\frac 12+i\vp_1^k(t)\right)\right|^2{\rm d}t,\ U\in\left(\left. 0,\frac{T}{\ln^2T}\right.\right] , \ T\to\infty ,
\end{split}
\ee
and, consequently, for $F(t)=1$ we have
\be \label{3.7}
\begin{split}
& g_{n+1}\frac{1}{U}\int_T^{T+U} \prod_{k=0}^n\left|\zeta\left(\frac 12+i\vp_1^k(t)\right)\right|^2{\rm d}t \sim \\
& \sim \prod_{l=1}^s g_l\frac{1}{U}\int_T^{T+U}\prod_{k=0}^{a_{j_l}-1}\left|\zeta\left(\frac 12+i\vp_1^k(t)\right)\right|^2{\rm d}t, \ T\to\infty .
\end{split}
\ee
\end{cor}

\subsection{}

Next, the set of all $p(n+1)-1$ asymptotic equalities (\ref{3.6}) generates the set of all
\bdis
\frac 12\left\{ p(n+1)-1\right\}\left\{ p(n+1)-2\right\},\ n\geq 2
\edis
asymptotic equalities between the right-hand sides of (\ref{3.6}).

\begin{mydef6}
Let us consider two partitions of the number 6:
\bdis
6=2+2+2,\ a_{j_l}=2;\quad 6=3+3,\ a_{j_l}=3 .
\edis
Then
\be \label{3.8}
\begin{split}
& \left\{g_2\frac{1}{U}\int_T^{T+U}\prod_{k=0}^1\left|\zeta\left(\frac 12+i\vp_1^k(t)\right)\right|^2{\rm d}t\right\}^3\sim \\
& \sim \left\{g_3\frac{1}{U}\int_T^{T+U}\prod_{k=0}^2\left|\zeta\left(\frac 12+i\vp_1^k(t)\right)\right|^2{\rm d}t\right\}^2 .
\end{split}
\ee
Since
\bdis
\begin{split}
& \int_T^{T+U}\prod_{k=0}^2\left|\zeta\left(\frac 12+i\vp_1^k(t)\right)\right|^2{\rm d}t= \\
& =\left|\zeta\left(\frac 12+i\vp_1^2(t_1)\right)\right|^2\int_T^{T+U}\prod_{k=0}^1\left|\zeta\left(\frac 12+i\vp_1^k(t)\right)\right|^2{\rm d}t,\ t_1\in (T,T+U)
\end{split}
\edis
then it follows from (\ref{3.8}) that
\be \label{3.9}
\begin{split}
& \left|\zeta\left(\frac 12+i\vp_1^2(t_1)\right)\right|^4\sim\frac{(g_2)^3}{(g_3)^2}\frac{1}{U}\int_T^{T+U}\prod_{k=0}^1\left|\zeta\left(\frac 12+i\vp_1^k(t)\right)\right|^2{\rm d}t= \\
& = \frac{(g_2)^3}{(g_3)^2}\left|\zeta\left(\frac 12+i\rho\right)\right|^2\left|\zeta\left(\frac 12+i\vp_1(\rho)\right)\right|^2,\ \rho\in (T,T+U) .
\end{split}
\ee
Hence, denoting
\bdis
\rho=\tau_0(T),\ \vp_1(\rho)=\tau_1(T),\ \vp_1^2(t_1)=\tau_2(T)
\edis
we obtain from (\ref{3.9}) the following corollary.
\end{mydef6}

\begin{cor}
There are values (a continuum set of these if $T\to\infty$)
\bdis
\tau_k(T)\in \left(\vp_1^k(T),\vp_1^k(T+U)\right),\ k=0,1,2
\edis
such that the asymptotic equality
\be \label{3.10}
\left|\zeta\left(\frac 12+i\tau_2\right)\right|^2\sim\frac{(g_2)^{3/2}}{g_3}\left|\zeta\left(\frac 12+i\tau_1\right)\right|\left|\zeta\left(\frac 12+i\tau_0\right)\right|,\ T\to\infty
\ee
holds true.
\end{cor}

\begin{remark}
The formula (\ref{3.10}) gives us a new kind of result about the distribution of the values
\bdis
\left|\zf\right|
\edis
with respect to the disconnected set
\bdis
\bigcup_{k=0}^2 \left[\vp_1^k(T),\vp_1^k(T+U)\right] .
\edis
\end{remark}

\begin{remark}
We hope that this example gives the sufficient information about the construction of the analogue of (\ref{3.10}) for the elements of the second set.
\end{remark}

\subsection{}

\begin{remark}
The formula (\ref{3.7}) gives us, in the direction from the left to the right, the factorization of the energy and, from the left to the right, we have the multiplicative synthesis of the
elementary energies in the following sense. For every fixed natural number $n+1$ we have the reservoir
\bdis
R=\bigcup_{L=1}^n \underbrace{\{ J_L,J_L,\dots,J_L\}}_{n_L}
\edis
where
\bdis
J_L=g_L\frac{1}{U}\int_T^{T+U}\prod_{k=0}^{L-1}\left|\zeta\left(\frac 12+i\vp_1^k(t)\right)\right|^2{\rm d}t,\ n_L=\left[\frac{n+1}{L}\right],
\edis
of the weighted elementary energies (of the same integral form). Now, the proper partition (\ref{3.2}) chooses
(an analogue of the \emph{shem ha-meforash of a golem}) corresponding elementary energies from the reservoir $R$, and by multiplication of these, we obtain
the resulting energy (of the same integral form).
\end{remark}

\section{Other factorizations of the integral (\ref{1.12})}

\subsection{}

We obtain by the substitution
\bdis
T\rightarrow \vp_1^k(T),\quad T+U\rightarrow \vp_1^k(T+U)
\edis
in the formula (\ref{1.7}) (see (\ref{2.3}), (\ref{2.4}))
\be \label{4.1}
\begin{split}
& \frac{1}{\vp_1^k(T+U)-\vp_1^k(T)}\int_{\vp_1^k(T)}^{\vp_1^k(T+U)}\left|\zf\right|^2{\rm d}t \sim \\
& \sim \frac{\vp_1^{k+1}(T+U)-\vp_1^{k+1}(T)}{\vp_1^k(T+U)-\vp_1^k(T)}\ln T,\ k=0,1,\dots \ ,n \\
& \ln\vp_1^k(T)\sim\ln T ,
\end{split}
\ee
and the multiple of these factors gives us the following formula
\be \label{4.2}
\begin{split}
& \prod_{k=0}^n\frac{1}{\vp_1^k(T+U)-\vp_1^k(T)}\int_{\vp_1^k(T)}^{\vp_1^k(T+U)}\left|\zf\right|^2{\rm d}t\sim \\
& \sim \frac{1}{U}\left\{\vp_1^{n+1}(T+U)-\vp_1^{n+1}(T)\right\}\ln^{n+1}T .
\end{split}
\ee
Next, we have from (\ref{2.1})
\be \label{4.3}
\begin{split}
& \frac{1}{U}\int_T^{T+U} F[\vp_1^{n+1}(t)]\prod_{k=0}^n\left|\zeta\left(\frac 12+i\vp_1^k(t)\right)\right|^2{\rm d}t\sim \\
& \sim \frac{1}{\vp_1^{n+1}(T+U)-\vp_1^{n+1}(T)}\int_{\vp_1^{n+1}(T)}^{\vp_1^{n+1}(T+U)}F(t){\rm d}t \ \times \\
& \times \frac{1}{U}\left\{\vp_1^{n+1}(T+U)-\vp_1^{n+1}(T)\right\}\ln^{n+1}T .
\end{split}
\ee
Hence, we obtain by (\ref{4.2}), (\ref{4.3}) the following corollary.

\begin{cor}
\be \label{4.4}
\begin{split}
& \frac{1}{U}\int_T^{T+U} F[\vp_1^{n+1}(t)]\prod_{k=0}^n\left|\zeta\left(\frac 12+i\vp_1^k(t)\right)\right|^2{\rm d}t\sim \\
& \sim \frac{1}{\vp_1^{n+1}(T+U)-\vp_1^{n+1}(T)}\int_{\vp_1^{n+1}(T)}^{\vp_1^{n+1}(T+U)}F(t){\rm d}t \ \times \\
& \sim \prod_{k=0}^n\frac{1}{\vp_1^{k+1}(T+U)-\vp_1^{k+1}(T)}\int_{\vp_1^k(T)}^{\vp_1^k(T+U)}\left|\zf\right|^2{\rm d}t
\end{split}
\ee
(comp. with the formula (\ref{1.10})).
\end{cor}

\subsection{}

In addition to (\ref{4.4}) there are also other (degenerate) factorizations. Namely, it follows from (\ref{4.1}) that
\be \label{4.5}
\begin{split}
& \frac{1}{\vp_1^{l+1}(T+U)-\vp_1^{l+1}(T)}\int_{\vp_1^l(T)}^{\vp_1^l(T+U)}\left|\zf\right|^2{\rm d}t\sim \ln T , \\
& l=0,1,\dots,n;\quad T\to\infty .
\end{split}
\ee
Now, we obtain from (\ref{2.1}) by (\ref{4.5}) the following corollary.

\begin{cor}
\be \label{4.6}
\begin{split}
& \int_T^{T+U} F[\vp_1^{n+1}(t)]\prod_{k=0}^n\left|\zeta\left(\frac 12+i\vp_1^k(t)\right)\right|^2{\rm d}t\sim \\
& \int_{\vp_1^{n+1}(T)}^{\vp_1^{n+1}(T+U)}F(t){\rm d}t \times \\
& \times \left\{\frac{1}{\vp_1^{l+1}(T+U)-\vp_1^{l+1}(T)}\int_{\vp_1^l(T)}^{\vp_1^l(T+U)}\left|\zf\right|^2{\rm d}t\right\}^{n+1} , \\
& l=0,1,\dots,n,\quad T\to\infty .
\end{split}
\ee
\end{cor}

For example, we obtain in the case (\ref{4.6}) with $l=0$ the following formula
\be \label{4.7}
\begin{split}
& \frac{\int_T^{T+U} F[\vp_1^{n+1}(t)]\prod_{k=0}^n\left|\zeta\left(\frac 12+i\vp_1^k(t)\right)\right|^2{\rm d}t}{\int_{\vp_1^{n+1}(T)}^{\vp_1^{n+1}(T+U)}F(t){\rm d}t} \sim \\
& \sim \left\{\frac{1}{\vp_1(T+U)-\vp_1(T)}\int_{T}^{T+U}\left|\zf\right|^2{\rm d}t\right\}^{n+1},\ T\to\infty .
\end{split}
\ee

\begin{remark}
The formula (\ref{4.7}) shows a curious effect: its left-hand side depends on all iterations $\vp_1^0,\dots,\vp_1^{n+1}$ while its right-hand side depends only the two first iterations
$\vp_1^0$ and $\vp_1^1$, where
\bdis
\vp_1^0(t)=t,\ \vp_1^1(t)=\vp_1(t) .
\edis
\end{remark}

\section{Generalization of some Solberg's formulae and the nonlocal and nonlinear interactions of the functions $|\zf|$ and $\arg\zeta(\frac 12+it)$}

\subsection{}

First, let us remind the Selberg's formula (see \cite{5}, p. 126, Theorem 6)
\be \label{5.1}
\begin{split}
& \int_T^{T+U}\{ S(t)\}^{2l}{\rm d}t\sim \frac{(2l)!}{l!(2\pi)^{2l}}U(\ln\ln T)^l, \\
& U\in \left[ T^{1/2+\delta},T\right],\quad T\to\infty ,
\end{split}
\ee
where $l$ is arbitrary fixed natural number, and
\be \label{5.2}
S(t)=\frac{1}{\pi}\arg\zf ,
\ee
with the $\arg$ function defined as usually (comp. \cite{6}, p. 179). By making use of the formulae (see (\ref{2.3}), (\ref{2.7}))
\be \label{5.3}
\vp_1^{n+1}(T)\sim T,\ \vp_1^{n+1}(T+U)-\vp_1^{n+1}(T)\sim U
\ee
we obtain from (\ref{5.1}) that
\be \label{5.4}
\begin{split}
& \int_{\vp_1^{n+1}(T)}^{\vp_1^{n+1}(T+U)}\{ S(t)\}^{2l}{\rm d}t\sim \\
& \sim \frac{(2l)!}{l!(2\pi)^{2l}} \left\{ \vp_1^{n+1}(T+U)-\vp_1^{n+1}(T)\right\}\left\{\ln\vp_1^{n+1}(T)\right\}^l\sim \\
& \sim \frac{(2l)!}{l!(4\pi)^{2l}}U(\ln\ln T)^l,\ T\to\infty ,
\end{split}
\ee
and consequently, we obtain from our Theorem in the case
\bdis
F(t)=S(t)
\edis
(see (\ref{5.2})) the following corollary.

\begin{cor}
\be \label{5.5}
\begin{split}
& \frac{1}{U}\int_T^{T+U}\left\{\arg\zeta\left(\frac 12+i\vp_1^{n+1}(t)\right)\right\}^{2l}\prod_{k=0}^n\left|\zeta\left(\frac 12+i\vp_1^k(t)\right)\right|^2{\rm d}t\sim \\
& \sim \frac{(2l)!}{4^ll!}(\ln\ln T)^l\ln^{n+1}T,\ U\in\left[ T^{1/2+\delta},\frac{T}{\ln^2T}\right],\ T\to\infty ,
\end{split}
\ee
for every fixed $l\in\mbb{N},\ n\in\mbb{N}_n$.
\end{cor}

\begin{remark}
The formula (\ref{5.5}) in the form
\bdis
\begin{split}
& \frac{1}{U}\int_T^{T+U}\left|\zf\right|^2\left\{\arg\zeta\left(\frac 12+i\vp_1^{n+1}(t)\right)\right\}^{2l}\prod_{k=0}^n\left|\zeta\left(\frac 12+i\vp_1^k(t)\right)\right|^2{\rm d}t\sim \\
& \sim \frac{(2l)!}{4^ll!}(\ln\ln T)^l\ln^{n+1}T
\end{split}
\edis
is simultaneously the generalization of the corresponding part of the Hardy-Littlewood integral (\ref{1.2}), (comp. (\ref{1.13})).
\end{remark}

\subsection{}

Next, let us remind another Selberg's formula (see \cite{5}, p. 130, Theorem 7)
\be \label{5.6}
\int_T^{T+U}\{ S_1(t)\}^{2l}{\rm d}t\sim d_lU,\ U\in\left[ T^{1/2+\delta},T\right],\ T\to\infty
\ee
where
\bdis
S_1(T)=\int_0^T S(t){\rm d}t.
\edis
We obtain, similarly to (\ref{5.4}), that
\bdis
\int_{\vp_1^{n+1}(T)}^{\vp_1^{n+1}(T+U)}\{ S_1(t)\}^{2l}{\rm d}t\sim d_l\left\{\vp_1^{n+1}(T+U)-\vp_1^{n+1}(T)\right\}d_l U ,
\edis
and consequently we obtain from the Theorem ($F(t)=S_1(t)$) the following corollary.

\begin{cor}
\be \label{5.7}
\begin{split}
& \frac{1}{U}\int_T^{T+U}\left\{ S_1[\vp_1^{n+1}(t)]\right\}^{2l}\prod_{k=0}^n\left|\zeta\left(\frac 12+i\vp_1^k(t)\right)\right|^2{\rm d}t\sim \\
& \sim d_l\ln^{n+1}T,\ U\in\left[ T^{1/2+\delta},\frac{T}{\ln^2T}\right],\ T\to\infty
\end{split}
\ee
for every fixed $l\in\mbb{N},\ n\in\mbb{N}_0$.
\end{cor}

\begin{remark}
We obtain from (\ref{5.5}), (\ref{5.7}) in the case $l=1,\ n=0$ the following minimal formulae
\be \label{5.8}
\begin{split}
& \frac{1}{U}\int_T^{T+U}\left\{\arg\zeta\left(\frac 12+i\vp_1(t)\right)\right\}^{2}\left|\zeta\left(\frac 12+it\right)\right|^2{\rm d}t\sim \\
& \sim \frac 12\ln\ln T\ln T, \\
& \frac{1}{U}\int_T^{T+U}\left\{ S_1[\vp_1(t)]\right\}^{2}\left|\zf\right|^2{\rm d}t\sim d_1\ln T,\ T\to\infty .
\end{split}
\ee
The simplest form of the nonlocal and nonlinear interactions of the pairs of signals
\bdis
\left\{\arg\zf,\left|\zf\right|\right\},\ \left\{ S_1(t),\left|\zf\right|\right\}
\edis
on the disconnected set
\bdis
\left[\vp_1(T),\vp_1(T+U)\right]\bigcup [T,T+U] .
\edis
is expressed by the formulae (\ref{5.8}).
\end{remark}

\section{Proof of the Theorem}

We use our formula (see \cite{3}, (6.2))
\be \label{6.1}
t-\vp_1(t)\sim (1-c)\pi(t);\quad \pi(t)\sim \frac{t}{\ln t},
\ee
where $\pi(t)$ is the prime-counting function.

\begin{remark}
The fundamental geometric property of the set of Jacob's ladders is expressed by the formula (\ref{6.1}): the difference of the abscissa and the ordinate of the point
\bdis
[t,\vp_1(t)]
\edis
of every curve
\bdis
y=\vp_1(t)
\edis
is asymptotically equal to $(1-c)\pi(t)$.
\end{remark}

\subsection{}

We have (see (\ref{6.1}), $t\to\vp_1^k(t)$)
\be \label{6.2}
\begin{split}
& \vp_1^k(t)-\vp_1^{k+1}(t)\sim (1-c)\frac{\vp_1^k(t)}{\ln\vp_1^k(t)},\ k=0,1,\dots,n+1 , \\
& t\in [T,T+U],\ U\in\left(\left. 0,\frac{T}{\ln^2T}\right.\right];\ \vp_1^0(t)=t
\end{split}
\ee
(comp. (\ref{1.11})) for arbitrary fixed $n\in\mbb{N}_0$, and (see (\ref{6.1}), (\ref{6.2}))
\be \label{6.3}
\begin{split}
& t\sim \vp_1^1(t)\sim\vp_1^2(t)\sim\dots\sim\vp_1^{n+1}(t),\ T\to\infty , \\
& t>\vp_1^1(t)>\vp_1^2(t)>\dots>\vp_1^{n+1}(t) .
\end{split}
\ee
Next, we have (see (\ref{6.2}), (\ref{6.3}))
\be \label{6.4}
\vp_1^k(t)-\vp_1^{k+1}(t)\sim\frac{t}{\ln t},\ k=0,1,\dots,n ,
\ee
and consequently we obtain by addition of (\ref{6.4})
\be \label{6.5}
\begin{split}
& t-\vp_1^{n+1}(t)\sim (1-c)(n+1)\frac{t}{\ln t} , \\
& \vp_1^{n+1}(t)\sim\left\{ 1-\frac{(1-c)(n+1)}{\ln t}\right\}t,\ 0<1-c<1,  \\
& \vp_1^{n+1}(t)>\left( 1-\frac{\epsilon}{2}\right)\left\{ 1-\frac{(1-c)(n+1)}{\ln t}\right\}t> \\
& > (1-\epsilon)t\geq (1-\epsilon)T,\ t\in [T,T+U] ,
\end{split}
\ee
i. e. from (\ref{6.3}), (\ref{6.5}) the properties (\ref{2.3}) follow. Especially, the following holds true
\be \label{6.6}
(1-\epsilon)T<\vp_1^{n+1}(T)<T .
\ee

\subsection{}

First of all we have (see (\ref{6.1}), (\ref{6.2}) -- condition for $U$)
\bdis
\begin{split}
& T-\vp_1^1(T)=\{ 1+o(1)\}(1-c)\frac{T}{\ln T} , \\
& T+U-\vp_1^1(T+U)=\{ 1+o(1)\}(1-c)\frac{T}{\ln T} ,
\end{split}
\edis
then
\bdis
\begin{split}
& \left( 1-\frac{1}{(2n+5)^2}\right)(1-c)\frac{T}{\ln T}<T+U-\vp_1^1(T+U)< \\
& < \left( 1+\frac{1}{(2n+5)^2}\right)(1-c)\frac{T}{\ln T}, \\
& \left( 1-\frac{1}{(2n+5)^2}\right)(1-c)\frac{T}{\ln T}<T-\vp_1^1(T)< \\
& < \left( 1+\frac{1}{(2n+5)^2}\right)(1-c)\frac{T}{\ln T},\quad T\to\infty .
\end{split}
\edis

Next, we have ($0<1-c<1$)
\bdis
\begin{split}
& \left| T+U-\vp_1^1(T+U)-\left\{ T-\vp_1^1(T)\right\}\right|< \\
& < \frac{2}{(2n+5)^2}(1-c)\frac{T}{\ln T}<\frac{2}{(2n+5)^2}\frac{T}{\ln T} ,
\end{split}
\edis
i. e.
\be \label{6.7}
\vp_1^1(T+U)-\vp_1^1(T)-U< \frac{2}{(2n+5)^2}\frac{T}{\ln T},
\ee
and (see (\ref{6.2}) -- the condition for $U$)
\bdis
0<\vp_1^1(T+U)-\vp_1^1(T)<\frac{2}{(2n+5)^2}\frac{T}{\ln T}+U<\frac{3}{(2n+5)^2}\frac{T}{\ln T} .
\edis
Hence (see (\ref{2.1}))
\be \label{6.8}
U\leq \frac{T}{\ln^2T}\ \Rightarrow \ \vp_1^1(T+U)-\vp_1^1(T)<\frac{3}{(2n+5)^2}\frac{T}{\ln T} .
\ee
Similarly, from the formula (see (\ref{6.2}))
\bdis
\vp_1^1(t)-\vp_1^2(t)\sim (1-c)\frac{t}{\ln t},\ t\to\infty ,
\edis
we obtain (comp. (\ref{6.7}), (\ref{6.8}))
\bdis
\begin{split}
& \vp_1^2(T+U)-\vp_1^2(T)<\frac{2}{(2n+5)^2}\frac{T}{\ln T}+\vp_1^1(T+U)-\vp_1^1(T)< \\
& < \frac{5}{(2n+5)^2}\frac{T}{\ln T}.
\end{split}
\edis

Next, if the estimate (the function $\vp_1^k(t)$ is increasing)
\bdis
\vp_1^k(T+U)-\vp_1^k(T)<\frac{2k+1}{(2n+5)^2}\frac{T}{\ln T}
\edis
holds true then we obtain by a similar way that
\bdis
\begin{split}
& \vp_1^{k+1}(T+U)-\vp_1^{k+1}(T)<\frac{2}{(2n+5)^2}\frac{T}{\ln T}+\frac{2k+1}{(2n+5)^2}\frac{T}{\ln T}< \\
& < \frac{2(k+1)+1}{(2n+5)^2}\frac{T}{\ln T} .
\end{split}
\edis
Hence, the following estimates hold true
\be \label{6.9}
\vp_1^{k}(T+U)-\vp_1^k(T)<\frac{2k+1}{(2n+5)^2}\frac{T}{\ln T}<\frac{1}{2n+5}\frac{T}{\ln T},\ k=1,\dots,n+1 ,
\ee
i. e. we have the inequalities (\ref{2.4}).

\subsection{}

Next we have (see (\ref{6.4}))
\bdis
\vp_1^k(T)-\vp_1^{k+1}(T)>\left( 1-\frac{1}{(2n+5)^2}\right)(1-c)\frac{T}{\ln T},
\edis
i. e. (see (\ref{6.9}))
\bdis
\begin{split}
& \vp_1^k(T)-\vp_1^{k+1}(T)=\vp_1^k(T)-\vp_1^{k+1}(T+U)+\vp_1^{k+1}(T+U)-\vp_1^{k+1}(T)> \\
& >\left( 1-\frac{1}{(2n+5)^2}\right)(1-c)\frac{T}{\ln T}, \\
& \vp_1^k(T)-\vp_1^{k+1}(T+U)>\left( 1-\frac{1}{(2n+5)^2}\right)(1-c)\frac{T}{\ln T}- \\
& -\left\{\vp_1^{k+1}(T+U)-\vp_1^{k+1}(T)\right\}> \\
& > \left( 1-c-\frac{1}{(2n+5)^2}\right)\frac{T}{\ln T}-\frac{1}{2n+5}\frac{T}{\ln T} = \\
& = \left( 1-c-\frac{1}{2n+5}-\frac{1}{(2n+5)^2}\right)\frac{T}{\ln T}\geq (1-c-0.24)\frac{T}{\ln T} > \\
& > 0.18\times\frac{T}{\ln T},\ k=0,1,\dots,n,\ n\geq 0
\end{split}
\edis
($c<0.58 \ \Rightarrow \ 1-c>0.42$), i. e. (\ref{2.5}) holds true.

\subsection{}

We use the Hardy-Littlewood-Ingham formula
\be \label{6.10}
\int_T^{T+U} Z^2(t){\rm d}t\sim U\ln T,\ U\in\left[ T^{1/3+\epsilon},\frac{T}{\ln^2T}\right]
\ee
(here the exponent $\frac 13$ is called the Balasubramanian exponent) in what follows. Next, we use also our formula (\ref{1.7})
\be \label{6.11}
\int_T^{T+U} Z^2(t){\rm d}t\sim\left\{ \vp_1(T+U)-\vp_1(T)\right\}\ln T.
\ee

We compare the formulae (\ref{6.10}) and (\ref{6.11}) in order to obtain
\bdis
\vp_1^1(T+U)-\vp_1^1(T)\sim U.
\edis
Similarly, by such a comparison in the cases
\bdis
T\to \vp_1^1(T),\ T+U\to \vp_1^1(T+U),\ \dots
\edis
where (see (\ref{2.3}))
\bdis
\ln\vp_1^k(T)\sim\ln T
\edis
we obtain
\be \label{6.12}
\vp_1^k(T+U)-\vp_1^k(T)\sim U,\ k=1,\dots,n+1,
\ee
i. e. the formula (\ref{2.7}) holds true.

\subsection{}

We have by (\ref{6.4})
\bdis
\vp_1^k(T)-\vp_1^{k+1}(T)\sim (1-c)\frac{T}{\ln T} ,
\edis
and consequently (see (\ref{6.2}) -- the condition for $U$, (\ref{6.12}))
\bdis
\begin{split}
& \vp_1^k(T)-\vp_1^{k+1}(T+U)\sim (1-c)\frac{T}{\ln T}-\left\{\vp_1^{k+1}(T+U)-\vp_1^{k+1}(T)\right\}= \\
& = (1-c)\frac{T}{\ln T}+\mcal{O}\left(\frac{T}{\ln^2T}\right)\sim (1-c)\frac{T}{\ln T},
\end{split}
\edis
i. e. the formula (\ref{2.8}) holds true. The formula (\ref{2.9}) follows from (\ref{2.8}).

\subsection{}

Let us remind (see \cite{4}, (9.1), (9.2)) that
\be \label{6.13}
\tilde{Z}^2(t)=\frac{{\rm d}\vp_1(t)}{{\rm d}t},\ \vp_1(t)=\frac 12\vp(t),\ t\geq T_0[\vp_1]
\ee
where
\be \label{6.14}
\begin{split}
& \tilde{Z}^2(t)=\frac{Z^2(t)}{2\Phi^\prime_\vp [\vp(t)]}=\frac{\left|\zf\right|^2}{\left\{ 1+\mcal{O}\left(\frac{\ln\ln t}{\ln t}\right)\right\}\ln t} , \\
& t\in [T,T+U],\ U\in \left(\left. 0,\frac{T}{\ln T}\right.\right] .
\end{split}
\ee
If we use the formula (\ref{6.13}) for the iterations (\ref{1.11}) we obtain
\be \label{6.15}
\prod_{k=0}^n \tilde{Z}^2[\vp_1^k(t)]=\frac{{\rm d}\vp_1^1}{{\rm d}t}\frac{{\rm d}\vp_1^2}{{\rm d}\vp_1^1}\cdots\frac{{\rm d}\vp_1^{n+1}}{{\rm d}\vp_n^1}=\frac{{\rm d}\vp_1^{n+1}}{{\rm d}t}
\ee
by the rule of differentiation of the composite function. Consequently we have (see (\ref{1.12}), (\ref{6.15}))
\bdis
\begin{split}
& \int_T^{T+U}F[\vp_1^{n+1}(t)]\prod_{k=0}^n\tilde{Z}^2[\vp_1^k(t)]{\rm d}t= \\
& =\int_T^{T+U}F[\vp_1^{n+1}(t)]{\rm d}\vp_1^{n+1}(t)=\int_{\vp_1^{n+1}(T)}^{\vp_1^{n+1}(T+U)} F(t){\rm d}t,
\end{split}
\edis
i. e.
\be \label{6.16}
\int_T^{T+U} F[\vp_1^{n+1}(t)]\prod_{k=0}^n\tilde{Z}^2[\vp_1^k(t)]{\rm d}t=\int_{\vp_1^{n+1}(T)}^{\vp_1^{n+1}(T+U)} F(t){\rm d}t .
\ee
Because (comp. (\ref{6.6}))
\bdis
(1-\epsilon)T<\vp_1^{n+1}(T)<T+U
\edis
we have the following
\be \label{6.17}
T' \in \left(\vp_1^{n+1}(T),T+U\right) \ \Rightarrow \ \ln T'=\ln T+\mcal{O}(1) .
\ee
Next, if we use the mean-value theorem on the left-hand side of the formula (\ref{6.16}), we obtain (see (\ref{2.3}) -- the first equality, (\ref{6.3}), (\ref{6.14}), (\ref{6.17}))
\be \label{6.18}
\begin{split}
& \int_T^{T+U}F[\vp_1^{n+1}(t)]\prod_{k=0}^n\tilde{Z}^2[\vp_1^k(t)]{\rm d}t\sim \\
& \sim \frac{1}{\ln^{n+1}T}\int_T^{T+U}\prod_{k=0}^n\left|\zeta\left(\frac 12+i\vp_1^k(t)\right)\right|^2{\rm d}t .
\end{split}
\ee
Hence, from (\ref{6.16}) and (\ref{6.18}) the asymptotic formula (\ref{2.1}) follows.

\thanks{I would like to thank Michal Demetrian for helping me with the electronic version of this work.}

\end{document}